\documentclass[a4paper,oneside]{amsart}
\usepackage{amssymb}
\usepackage{amsfonts}
\usepackage{amsmath}

\newtheorem{theorem}{Theorem}[section]

\newtheorem{corollary}[theorem]{Corollary}

\newtheorem{lemma}[theorem]{Lemma}

\newtheorem{remark}[theorem]{Remark}

\title{Analytic functions with hyperbolic range and Bohr's inequality}
\author{Yusuf Abu Muhanna$^{1}$ and Issam Louhichi$^{2}$
\\
\\
$^{1,2}$ Department of Mathematics \& Statistics, College of Arts \& Sciences
\\
American University of Sharjah, P.O.Box 26666, Sharjah
\\
United Arab Emirates
\\
Email: ymuhanna@aus.edu, ilouhichi@aus.edu (corresponding author)}

\begin{document}
\maketitle

\begin{abstract}
	We use properties of the hyperbolic metric and properties of the modular
	function to show that the Bohr's radius for covering maps onto hyperbolic
	domains is $\geq e^{-\pi }.$ This includes almost all known classes of
	analytic functions.
\end{abstract}


\section{Introduction}

In this article, we shall investigate Bohr's phenomenon for spaces of
analytic functions into hyperbolic domains. Throughout, we let $U$ to be the
unit disk of the complex plane $\mathbb{C}$, $D$ a general domain and
 $f:U\rightarrow D$  an analytic function. We shall mainly focus on functions that map into hyperbolic domains.
A domain $D$, is called hyperbolic if its complement in $\mathbb{C}$ contains at least two finite points. As for example,  almost all known classes of
analytic functions contain maps mapping into hyperbolic domains. Finally,
we recall that an analytic function $\varphi :U\rightarrow U$ is called Schwarz if $
\varphi (0)=0.$

This manuscript is structured as follows: Subsection 1.1 contains some basic properties of the Bohr's operator. Subsection
1.2 contains the uniformization theorem, the hyperbolic metric and its basic properties. Subsection 1.3 is about the modular function and its basic properties. In addition to that, some basic consequences are listed without proofs. Those well-known results shall be used
 in Section 2. Section 2 contains the main result which is Theorem 2.1 and some direct
consequences. Section 3 contains some applications of Theorem 2.1 to harmonic maps.

\subsection{Bohr's operator}

We say that a class $\digamma $ of analytic functions  satisfies the Bohr's
phenomenon, if for a function $f(z)=\sum\limits_{0}^{\infty }a_{n}z^{n}$ in $%
\digamma $, the Bohr's operator at $z$ 
\begin{equation*}
M(f)=\sum\limits_{0}^{\infty }\left\vert a_{n}z^{n}\right\vert
\end{equation*}%
is uniformly bounded on some closed disk $\{|z|\leq \rho \}$, with $\rho >0$. The largest such radius $\rho $ is called the Bohr's radius for the class $%
\digamma$. See [1] and [5] for more details.

It is well-known that the Bohr's operator satisfies the following
properties:
\begin{itemize}
\item[i)] $M(f+g)\leq M(f)+M(g),$

\item[ii)] $M(fg)\leq M(f)M(g),$

\item[iii)] $M(1)=1.$
\end{itemize}
These properties make $\digamma $  a Banach algebra, with norm $M(f)$.
We recall the following two properties for a Banach algebra:
\begin{itemize}
\item[1)] \textit{Bohr's Theorem [4]}: If $|f(z)|<1,$ for all $z\in U,$ then 
$M(f)<1$, when $|z|<1/3.$
\item[2)] \textit{Von Neumann's inequality}: $|p(f(z))|<||p||_{\infty}$, where 
$p(z)$ is a polynomial.
\end{itemize}
In fact, von Neumann showed that the above inequality "2)" is true for the space of bounded operators on a Hilbert space (see [4, 6]). Later Dixon [6] showed that in the space 
\begin{equation*}
l_{\beta }^{1}=\left\{ x=(x_{1},x_{2},x_{3},\ldots ):\frac{1}{\beta }%
\sum_{j=1}^{\infty }|x_{j}|<\infty \right\} ,
\end{equation*}
the von Neumann's inequality is satisfied for $0<\beta <1/3$ but not for $\beta \geq 1/3.$

From now on and through out the rest of the paper, we deal with analytic functions $f$ on $U$ that are missing two points. Those functions are usually called the hyperbolic functions.

\subsection{Hyperbolic metric.}
The following theorem is central for the study of hyperbolic metrics on
hyperbolic domains.
\begin{theorem}
\bigskip (Uniformization Theorem [5]) If $D$ is hyperbolic, then there is a universal cover $F$ (conformal) from $U$ onto $D.$ This cover is unique with the normalization $F(0)=a$ and $F^{\prime }(0)>0,$ for some $a\in D.$
\end{theorem}

\begin{corollary} ([5]) If $f(z)$ is analytic and maps $U$ into $D$, then there is a Schwarz function $\varphi (z)$ so that $f(z)=F(\varphi (z)),$ where $F$ is a covering map of $D.$
\end{corollary}
The hyperbolic metric on $U$ (as mentioned in [5]) is given by 
\begin{equation*}
\lambda _{U}(z)=\frac{1}{1-|z|^{2}}.
\end{equation*}

If we denote the universal covering of a hyperbolic domain $D$ by $F,$ then $%
F$ generates a hyperbolic metric on $D$ defined by (see [5, p.43] and [9]) 
\begin{equation}
\lambda _{D}(F(z))=\frac{1}{F^{\prime }(z)}\frac{1}{1-|z|^{2}}.
\end{equation}%
We recall that if $d(w,\partial D)$ is the distance from $w$ to the boundary of 
$D$ and if $D$ is a hyperbolic domain, then it is known that (see [5]) 
\begin{equation}
d(w,\partial D)\lambda _{D}(w)\leq 1,
\end{equation}%
and if in addition $D$ is simply connected, we have that (see [5]) 
\begin{equation*}
\frac{1}{4}\leq d(w,\partial D)\lambda _{D}(w)\leq 1.
\end{equation*}%
Thus, for a covering map $F(z)$, (1) and (2) imply that 
\begin{equation}
d(F(z),\partial D)\leq F^{\prime }(z)(1-|z|^{2}),
\end{equation}%
and when $F(z)$ is univalent 
\begin{equation}
\frac{1}{4}F^{\prime }(z)(1-|z|^{2})\leq d(F(z),\partial D)\leq F^{\prime
}(z)(1-|z|^{2}).
\end{equation}

In the following theorems by David Minda [9], "conformal" means non-vanishing derivative.

\begin{theorem}
Let $D$ be a hyperbolic domain, with hyperbolic metric $\lambda (w),$ and
let $f:U\rightarrow D$ be a conformal map that is onto. Then 
\begin{equation*}
\lambda (f(z))=\frac{1}{|f^{\prime }(z)|}\frac{1}{1-|z|^{2}}.
\end{equation*}
\end{theorem}

\begin{theorem}
Let $f(z)=\sum\limits_{0}^{\infty }a_{n}z^{n}$ be analytic on $U$ and
suppose that $f(U)=D$ misses at least two finite points $a,b$ with $b\neq a$
(hyperbolic), then 
\begin{eqnarray*}
h(z) &=&z\frac{f(z)-a}{b-a} \\
&=&z\frac{a_{0}-a}{b-a}+\sum_{1}^{\infty }\frac{a_{n}z^{n+1}}{b-a}, \\
\end{eqnarray*}
is 0 only at 0.
\end{theorem}

\subsection{The Modular function and consequences}
The modular function is defined by 
\begin{equation*}
J(z)=16z\prod\limits_{1}^{\infty }\left[ \frac{1+z^{2n}}{1+z^{2n-1}}\right]
^{8},
\end{equation*}%
where $J(z)$ is $0$ only at $0$ and $J\neq 0,1,\infty $ on $\{|z|>0\}$. Note
that 
\begin{equation}
-J(-z)=z\sum\limits_{0}^{\infty }M_{n}z^{n},\text{ with }M_{n}>0\text{ for
all }n.
\end{equation}
For more details on the modular function, see $[11,12]$. This function is somehow similar to the Koebe function for a large function. We might call it
the "large Koebe". Now, from (5) we can immediately deduce that 
\begin{equation}
\underset{|z|=r}{\max }|J(z)|=|J(-r)|.
\end{equation}

\begin{lemma} ([11]) If $h(z)=\sum\limits_{0}^{\infty }a_{k}z^{k}$ is only $0$ at $0$ and if $a$ is in the complement of
$h(U)$, then $h(z)/a$ is subordinate to $-J(-z)$ and  $|a_{k}|\leq 16|a|M_{k},$
for all $k$.
\end{lemma}

We shall also use the following results. For the sake of completeness, we choose to list them without proofs. Proofs can be found in [11].
\begin{lemma}
The modular function $J$ satisfies the following properties:
\begin{itemize}
\item[a)] [11, p.85] $J(z)$ has radius of univalence $e^{-\pi /2},$\\
and
\item[b)] 
$|J(-e^{-\pi })|=1,\text{ }|J(e^{-\pi })|=1/2,\text{ and }\underset{|z|\leq
e^{-\pi }}{\max }|J(z)|=1.$
\end{itemize}
\end{lemma}
The above lemma and the fact that $e^{-\pi /2}>e^{-\pi }$ imply the following well-known results.

\begin{lemma}
\lbrack 11, Lemma 1, p.83] If $h(z)$ is analytic on $U$ and bounded with $h(0)=0$
and if $h(z)\neq 0$ for all $0<|z|<1$, then $h(z)$ is univalent in $|z|<\rho
,$ where 
\begin{equation*}
\rho =1+\alpha -\sqrt{\left( 1+\alpha \right) ^{2}-1}<0.26
\end{equation*}%
and $e^{-\alpha }=\left( \frac{h(z)}{z}\right) (0).$
\end{lemma}

\begin{corollary}
Let $h(z)$ be as in Lemma 1.7. Assume that $|z|\leq e^{-\pi }$. Then $h(z)$ is univalent in $|z|<e^{-\pi }.$
\end{corollary}

\section{Main results}

Below is our main theorem.
\begin{theorem}
Let $f(z)=\sum\limits_{0}^{\infty }a_{n}z^{n}$ be analytic on $U$ and
suppose that $f(U)=D$ misses at least two finite points $a\neq b$ (hyperbolic).
Then 
\begin{equation*}
\sum\limits_{1}^{\infty }|a_{n}z^{n}|\leq 2d(f(0),\partial f(U)),
\end{equation*}%
for $|z|<e^{-\pi }.$
\end{theorem}

\begin{proof}
$f(z)$ is subordinate to some covering map $F$ from $U$ onto $D$. The
existence of such covering is ensured by Theorem 1.1. Without loss of
generality, we may assume that 
\begin{equation*}
|f(0)-a|=|a-a_{0}|=d(f(0),\partial f(U)).
\end{equation*}%
For otherwise, replace $a$ by the nearest point in the complement of $D$ to $f(0)$. As in Theorem 1.4, the function $h$ defined by 
\begin{equation}
h(z)=z\frac{f(z)-a}{b-a},
\end{equation}%
is $0$ only at $0$. From (7), we can write that
\begin{eqnarray*}
f(z) &=&(b-a)\frac{h(z)}{z}+a \\
&=&f_{1}(z)+f(0),
\end{eqnarray*}%
where $f_{1}(z)=\sum\limits_{1}^{\infty }a_{n}z^{n}$. So $f_1(z)=f(z)-f(0)=(b-a)\frac{h(z)}{z}+a-f(0)$. Thus, using
the properties of the Bohr's operator, we obtain that 
\begin{equation}
M(f_{1})\leq |b-a|M\left(\frac{h(z)}{z}\right)+|a-f(0)|.
\end{equation}
Again, using the properties of the Bohr's operator, we deduce  from (7) and (8) that 
\begin{equation*}
M(zf)\leq |b-a|M(h)+M(az).
\end{equation*}%
Next, we shall estimate $M(h)$. Since $f(z)=\sum\limits_{0}^{\infty
}a_{n}z^{n}$, we have 
\begin{eqnarray}
h(z) &=&\frac{a_{0}-a}{a-b}z+\frac{z}{b-a}\sum\limits_{1}^{\infty }a_{n}z^{n}
\\
&=&\sum\limits_{1}^{\infty }c_{n}z^{n}  \notag
\end{eqnarray}%
with $|c_{n}|=|\frac{a_{n-1}}{b-a}|\text{ for }n>1,\text{ and }|c_{1}|=|%
\frac{a_{0}-a}{b-a}|.$ If we denote by $\delta =d(0,\partial h(U))>0$ and
because $\frac{h(z)}{\delta}$ is $0$ only at $0$ and misses $1$, we have that
\begin{equation}
\frac{h}{\delta} =-J(-\omega ),\text{ with }\omega \text{ being  Schwarz}.
\end{equation}%
Since the coefficients of $-J(-z)$ are convex increasing,
Lemma 1.5 yields 
\begin{equation*}
|c_{n}|=\left\vert \frac{a_{n-1}}{b-a}\right\vert <\delta M_{n}.
\end{equation*}%
On the other hand, by Corollary 1.8, the function $h_{1}(z)=h(e^{-\pi }z)$ is
univalent and bounded in $U$. Then from (9), we have 
\begin{eqnarray}
h_{1}(z) &=&\frac{a_{0}-a}{a-b}e^{-\pi }z+\frac{z}{b-a}\sum\limits_{1}^{%
\infty }a_{n}e^{-n\pi }z^{n}  \notag \\
&=&\sum\limits_{1}^{\infty }c_{n}e^{-n\pi }z^{n}.
\end{eqnarray}%
Let $\delta _{1}=d(0,h_{1}^{C}(U)),$ where $h_{1}^{C}(U)$ is the complement
of $h_{1}(U)$. Then $\frac{h_{1}(z)}{\delta _{1}}$ is subordinate to $-J(-z)$. Moreover,
for $|z|=1$, Lemma 1.6 and (11) imply 
\begin{equation}
M(h_{1})=\sum\limits_{1}\left\vert c_{n}\right\vert e^{-\pi (n)}\leq \delta
_{1}(-J(-e^{-\pi }))\leq \delta _{1},
\end{equation}%
and for $|z|=e^{-\pi },$ we have 
\begin{eqnarray}
M(zf(z)) &=&|a_{0}|e^{-\pi }+\sum\limits_{2}|a_{n-1}|e^{-n\pi } \\
&\leq &|b-a|M(h)+M(az)  \notag \\
&=&|b-a|M(h_{1})+M(az).  \notag \\
&\leq &\delta _{1}|b-a|+|a|e^{-\pi }.  \notag
\end{eqnarray}%
Thus 
\begin{equation*}
|a_{0}|+\sum\limits_{2}|a_{n-1}|e^{-(n-1)\pi }\leq e^{\pi }\delta
_{1}|b-a|+|a|,
\end{equation*}%
and so 
\begin{eqnarray*}
\sum\limits_{2}|a_{n-1}|e^{-(n-1)\pi } &\leq &e^{\pi }\delta
_{1}|b-a|+|a|-|a_{0}| \\
&\leq &e^{\pi }\delta _{1}|b-a|+|a-a_{0}|,
\end{eqnarray*}%
or when replacing $e^{-\pi }$ by $|z|$ 
\begin{eqnarray}
\sum\limits_{2}^{\infty }|a_{n-1}||z|^{n-1} &\leq &|a-a_{0}|\left( \delta
_{1}\frac{|b-a|}{e^{-\pi }|a-a_{0}|}+1\right)   \notag \\
&\leq &|a-a_{0}|\left( \frac{\delta _{1}}{|h_{1}^{\prime }(0)|}+1\right),
\end{eqnarray}%
where the last inequality is obtained because $$h'_1(z)=\frac{a_0-a}{a-b}e^{-\pi}+\frac{1}{b-a}\sum_{n\geq 1}(n+1)a_ne^{-n\pi}z^n,$$ and so $h'_1(0)=\frac{a_0-a}{a-b}e^{-\pi}$.
Note that $h^{\prime }(0)=\frac{a_{0}-a}{b-a}=J(e^{-\alpha })<J(-e^{-1})$.
We recall from (12), that $\delta _{1}J(e^{-\pi }\psi (z))=h_{1}(e^{-\pi }z),
$ where $\psi (z)$ is Schwarz. As $\delta _{1}J(e^{-\pi }z)$ and $h(e^{-\pi
}z)$ are $0$ only at $z=0,$ so is $\psi $. As $\rho $ in Lemma 1.7 is
greater than $e^{-\pi },$ $\psi (e^{-\pi }z)$ is univalent and consequently, 
$h(e^{-\pi }z)=h_{1}(z)$ is univalent. Hence, by (4) and noting that $%
\lambda _{h_{1}(U)}(z)=\left\vert e^{-\pi }h^{\prime }(e^{-\pi
}z)\right\vert (1-|z|^{2})$ and $d(h_{1}(z),\partial (h_{1}(z)))=\delta _{1},
$ we deduce that 
\begin{equation*}
\frac{1}{4}\left\vert e^{-\pi }h^{\prime }(e^{-\pi }z)\right\vert
(1-|z|^{2})\leq d(h_{1}(z),\partial (h_{1}(z)))\leq \left\vert e^{-\pi
}h^{\prime }(e^{-\pi }z)\right\vert (1-|z|^{2}).
\end{equation*}
Thus, at $z=0,$ we obtain 
\begin{equation*}
\frac{1}{4}\left\vert e^{-\pi }h^{\prime }(0)\right\vert \leq
d(h_{1}(0),\partial (h_{1}(z)))\leq \left\vert e^{-\pi }h^{\prime
}(0)\right\vert,
\end{equation*}
and%
\begin{equation*}
\frac{\delta _{1}}{e^{-\pi }|h^{\prime }(0)|}=\frac{\delta _{1}}{%
|h_{1}^{\prime }(0)|}\leq 1,
\end{equation*}
where the last inequality is due to the fact that $|h'_1(0)|=e^{-\pi}|h'(0)|$.
Consequently, $(14)$ becomes%
\begin{eqnarray*}
\sum\limits_{2}^{\infty }|a_{n-1}||z|^{n-2} &\leq &|a-a_{0}|(1+1) \\
&= &2|a-a_{0}| \\
&=&2 d(f(0),\partial f(U)).
\end{eqnarray*}

Therefore, since $\frac{f(z)-a_{0}}{z}=\sum\limits_{2}^{\infty
}a_{n-1}z^{n-2}, $ we conclude that 
\begin{equation*}
M(f(z)-a_{0})\leq 2d(0,\partial f(U)).
\end{equation*}
\end{proof}

As a direct consequence of our main result, we have the following straightforward corollaries for which we omit the proofs.

\begin{corollary}
If $\digamma $ is a class of uniformly bounded analytic functions on $U$, then
all functions in the class $\digamma $ miss same two points, and
hence $\digamma $ has Bohr's radius $\geq e^{-\pi }.$

\end{corollary}


\begin{remark}
\begin{itemize}

\item[1)] The authors strongly believe that the constant $2$ in the
inequality of the main result is not sharp and shall be reduced to $1$.

\item[2)] Let $f$ and $h$ be as in the proof of our main result. Then 
\begin{eqnarray*}
d(f(0),\partial f(U)) &\leq &|f^{\prime }(0)| \\
&=&|b-a||\left(\frac{h}{z}\right)^{\prime }(0)| \\
&=&|b-a||J^{\prime }(-e^{-i\beta }e^{-\alpha })e^{-\alpha }2\alpha |.
\end{eqnarray*}%
Hence 
\begin{equation*}
d(f(0),\partial f(U))\leq |b-a|J^{\prime }(-e^{-\alpha })e^{-\alpha }2\alpha.
\end{equation*}
\end{itemize}
\end{remark}
The following direct corollary shall be used in the next section.

\begin{corollary}
If $f(z)=\sum\limits_{0}^{\infty }a_{n}z^{n}$ is analytic on $U$ and if $%
f(U) $ misses at least two points, with $d(f(0),\partial f(U))<1$, then $%
f(z) $ satisfies the von Neumann's inequality for $|z|\leq e^{-\pi
}/3=1.4405\times 10^{-2}$. In other words, for any polynomial $p(z)$, we
have 
\begin{equation*}
p(f(z))\leq ||p||_{\infty }.
\end{equation*}
\end{corollary}

\section{Harmonic maps}

In this section $M$ denotes the Bohr's operator as mentioned previously in the introduction. We recall that for harmonic function $f=f_1+\overline{f_2}$ on $U$ (i.e., $f_1$ and $\overline{f_2}$ are analytic on $U$), $M(f)$ is defined to be $M(f)=M(f_1)+M(\overline{f_2})$. The following theorem is a consequence of Theorem 2.1.
\begin{theorem}
Let $f(z)=h(z)+\overline{g(z)}$ be a harmonic function on $U,$ where $h(z)$
and $g(z)$ are analytic. Assume that $h(z)$ maps $U$ onto a hyperbolic domain, $h(0)=a_{0}$, $g(0)=0$%
, and $g^{\prime }(z)=\mu (z)h^{\prime }(z),$ where $\mu (z)$, Schwarz, is
the dilatation of the map $f(z)$. Then 
\begin{equation*}
M(f-a_{0})\leq 4d(a_{0},\partial (h(U))),
\end{equation*}%
for $|z|\leq e^{-\pi }/3=1.4405\times 10^{-2}.$
\end{theorem}

\begin{proof}
Using the properties of the operator $M$, mentioned in Section 1.1, we have 
\begin{equation*}
M(g)=\int\limits_{0}^{r}M(g^{\prime })dr\leq \int\limits_{0}^{r}M(\mu
)M(h^{\prime })dr.
\end{equation*}%
Thus for $|z|<1/3$, $M(\mu )<1$ and so
\begin{equation*}
M(g)\leq \int\limits_{0}^{r}M(h^{\prime })dr=M(h)-|a_{0}|=M(h-a_{0}).
\end{equation*}%
In particular, for $|z|<e^{-\pi },$ Theorem 2.1 gives 
\begin{equation*}
M(g)<2d(a_{0},\partial h(U)),
\end{equation*}%
and therefore 
\begin{eqnarray*}
M(f-a_{0}) &=&M(h-a_{0})+M(g) \\
&\leq &4d(a_{0},\partial (h(U)).
\end{eqnarray*}
\end{proof}

\section{No conflict of interest} On behalf of all authors, the corresponding author states that there is no conflict of interest. 
\section{Compliance with Ethical Standards} Not applicable. 
\section{Data availability statement} The authors did not analyze or generate any data sets, because this work proceeds within a theoretical and mathematical approach.

\end{document}